\documentclass[dvipdfmx]{article}

\usepackage{amsmath,amsfonts,amssymb,amsthm,pb-diagram,picinpar,graphicx,color}
\usepackage{algorithm}
\usepackage{algorithmic}
\usepackage{colortbl}
\usepackage{enumerate}
\usepackage{hyperref}
\usepackage{amscd}
\usepackage{bm}
\usepackage{tikz}
\usepackage{pgf}
\usepackage[T1]{fontenc}
\usepackage{enumitem}
\usetikzlibrary{shapes,matrix,calc,arrows,decorations.markings,knots,patterns,intersections,cd}

\newcommand{\CC}{\mathbb{C}}
\newcommand{\QQ}{\mathbb{Q}}
\newcommand{\ZZ}{\mathbb{Z}}
\newcommand{\PP}{\mathbb{P}}
\newcommand{\FF}{\mathbb{F}}

\newcommand{\mcC}{\mathcal{C}}   
\newcommand{\mcB}{\mathcal{B}}

\newcommand{\mcL}{\mathcal{L}}

\newcommand{\mcO}{\mathcal{O}}
\newcommand{\mcP}{\mathcal{P}}
\newcommand{\mcQ}{\mathcal{Q}}
\newcommand{\mcT}{\mathcal{T}}

\newcommand{\fd}{\mathfrak{d}}

\newcommand{\vc}{\boldsymbol {c}}

\newcommand{\I}{\mathop{\mathrm{I}}\nolimits}
\newcommand{\III}{\mathop{\mathrm{III}}\nolimits}

\newcommand{\pic}{\mathop{\mathrm{Pic}}\nolimits}

\newcommand{\Div}{\mathop{\mathrm{Div}}\nolimits}

\newcommand{\Gal}{\mathop{\mathrm{Gal}}\nolimits}
\newcommand{\MW}{\mathop{\mathrm{MW}}\nolimits}
\newcommand{\NS}{\mathop{\mathrm{NS}}\nolimits}

\newcommand{\Supp}{\mathop{\mathrm{Supp}}\nolimits}

\newcommand{\sr}{\mathop{\mathrm{sr}}\nolimits}
\newcommand{\sing}{\mathop{\mathrm{Sing}}\nolimits}
\newcommand{\Red}{\mathop{\mathrm{Red}}\nolimits}

\newcommand{\ord}{\mathop{\mathrm{ord}}\nolimits}

\newcommand{\Kbar}{\overline{K}}

\newtheorem{thm}{Theorem}[section]
\newtheorem{thm0}{Theorem}
\newtheorem{lem}[thm]{Lemma}     
\newtheorem{cor}[thm]{Corollary}
\newtheorem{prop}[thm]{Proposition}
\newtheorem{prop0}{Proposition}

\theoremstyle{definition}
\newtheorem{defin}[thm]{Definition}
\newtheorem{defin0}[prop0]{Definition}
\newtheorem{ex}[thm]{Example}

\theoremstyle{remark}
\newtheorem{rem}[thm]{Remark}

\begin{document}
\title{
Trisections on Certain Rational Elliptic Surfaces and Families of Zariski Pairs Degenerating to the same Conic-line Arrangement}
\date{\today}

\author{S. BANNAI\footnote{Partially supported by Grant-in-Aid for Scientific Research C (18K03263)}
, N. KAWANA, R. MASUYA and  H. TOKUNAGA\footnote{Partially supported by Grant-in-Aid for Scientific Research C (20K03561)}}

\maketitle

\begin{abstract}
In this paper, we study the geometry of trisections on certain rational elliptic surfaces. We utilize   Mumford representations of semi-reduced divisors in order to construct trisections and related plane curves with interesting properties explicitly. As a result we are able to  construct new examples of Zariski pairs. Especially, we show the existence of a family of Zariski pairs that degenerate to the same conic-line arrangement.
\end{abstract}

\section*{Introduction}

Let $\varphi: S \to C$ be an elliptic surface over a smooth projective curve $C$ satisfying the conditions as follows:
\begin{enumerate}
\item[(i)] The morphism $\varphi$ is relatively minimal.
\item[(ii)]  There exists at least one section $O : C \to S$. We identify $O$ with its image, i.e.,  a curve on $S$ meeting
any fiber of $\varphi$ at  one point.

\item[(iii)] There exists  at least one singular fiber.
\end{enumerate}
Under these assumptions, the N\'eron-Severi group $\NS(S)$ of $S$ is finitely generated and
torsion-free by \cite[Theorem 1.2]{shioda90}.  Let $\MW(S)$ be the set of sections of $\varphi
: S \to C$. As $O \in \MW(S)$, $\MW(S) \neq \emptyset$. Note that for a section $s: C \to S$, we identify $s$ with
its image on $S$.

Let $E_S$ be the generic fiber of $\varphi : S \to C$. $E_S$ can be regarded as a curve of 
genus $1$ over the rational function field, $\CC(C)$, of $C$. Under our setting, $S$ is known as
the Kodaira-N\'eron model of $E_S$ and we can identify 
the set of $\CC(C)$-rational points $E_S(\CC(C))$ of $E_S$  with $\MW(S)$. Here we define the 
zero of $E_S(\CC(C))$ to be the $\CC(C)$-rational point given by the restriction $O$ to $E_S$. We also
denote it by $O$ for simplicity. For $P \in E_S(\CC(C))$, we denote the corresponding section by
$s_P$.

 Let $D$ be a divisor on $S$ and let $\fd_D$ be the restriction of $D$ to $E_S$. $\fd_D$ is a
divisor on  $E_S$ defined over $\CC(C)$. By Abel's theorem, we have $P_D \in E_S(\CC(C))$
from $\fd_D$ and the corresponding section $s_{P_D}$ which we also denote by $s(D)$. This relates $\NS(S)$ 
with $E_S(\CC(C))$ and as  is shown in \cite{shioda90}, we have

\begin{thm0}\label{thm:shioda}{{\rm (\cite[Theorem 1.3]{shioda90})}
\[
\bar{\psi}:\NS(S)/T_{\varphi} \cong E_S(\CC(C)),
\]
where $T_{\varphi}$ is a subgroup of $\NS(S)$ generated by $O$ and irreducible curves contained in fibers.
(A more detailed description for $T_{\varphi}$ can be found in Section~\ref{sec:elliptic-surfaces}.)
}
\end{thm0}

In \cite[Lemma 8.2]{shioda90}, Shioda defined a homomorphism $\phi: E_S(\CC(C)) \to
\NS(S)\otimes\QQ$ by which he defined a pairing $\langle \bullet, \bullet\rangle$ on  $E_S(\CC(C))$ called the {\it height pairing} and went on to define a   lattice structure on the free part of $E_S(\CC(C))$, namely the Mordell-Weil
lattice. In \cite[Section~2]{bannai-tokunaga20}, the first and the last author defined a homomorphism
$\phi_o: \Div(S) \to \NS(S)\otimes\QQ$ such that $\phi(P_D) = \phi_o(D)$. With $\phi, \phi_o$ and 
Theorem~\ref{thm:shioda}, the geometric relation between $D$ and $P_D$ can be investigated.

In this article, we consider curves whose classes are mapped to a given $P \in E_S(\CC(C))$ in the case when $S$ is rational, in which case $C = \PP^1$. 
We need to introduce some terminologies, before we go on to explain our problem.
 
 For a divisor $D$ on $S$,  we say that $D$ is {\it vertical} if all irreducible components are contained in 
fibers. $D$ is said to be {\it horizontal} if $D$ contains no vertical components.  

\begin{defin0}\label{def:m-section}{\rm A horizontal  divisor
$D$ is said to be an $m$-section if $D$ meets in  $m$ distinct points with a  general fiber. Furthermore if $D\cdot O=0$, $D$ is called an {\it integral} $m$-section.}
\end{defin0}

Note that an $1$-section
is just an ordinary  section. We call $D$ a {\it bisection} or {\it trisection} if $m = 2$ or $=3$, respectively. 

The purpose of this article is to study trisections of $S$ to find plane curves in which we
are interested, in particular curves that are interesting from the viewpoint of  the embedded topology of plane curve arrangements with irreducible components of  low degree.
The first and the last authors have been investigating bisections and trisections along this line when 
$S$ is  a certain rational elliptic surface associated to
a reduced plane quartic $\mcQ$ and a smooth point $z_o$.  In order to make our problem clear we explain
how we obtain such a surface more precisely.

Let $\mcQ$ be a reduced quartic which is not the union of  $4$ concurrent lines. Let $z_o$ be a smooth point of $\mcQ$.
We associate a rational elliptic 
surface, $S_{\mcQ, z_o}$, to $\mcQ$ and $z_o$ as follows: 
\begin{enumerate}
\item[(i)]  Let $f'_{\mcQ}: S'_{\mcQ} \to \PP^2$ be the double cover branched along $\mcQ$. We denote its 
canonical resolution by $\mu : S_{\mcQ} \to
S'_{\mcQ}$ (see \cite{horikawa} for the canonical resolution). Note that $S_{\mcQ}$ is a rational surface.

\item[(ii)] Let $\Lambda_{z_o}$ be the pencil of curves of genus $1$ on $S_{\mcQ}$ induced from the pencil of 
lines through $z_o$. The pencil $\Lambda_{z_o}$ has a unique  base point $f_{\mcQ}^{-1}(z_o)$ of multiplicity $2$.
Let $\nu_{z_o}: S_{\mcQ, z_o} \to S_{\mcQ}$ be the resolution of the indeterminacy for the rational map
given by $\Lambda_{z_o}$. The induced morphism
 $\varphi_{\mcQ, z_o} : S_{\mcQ, z_o} \to \PP^1$ is an elliptic fibration.  $\varphi_{\mcQ, z_o}$ 
 satisfies the three conditions at the beginning of the Introduction. Thus we have a rational elliptic surface
 $S_{\mcQ, z_o}$ and the diagram below:
  \[
\begin{CD}
S'_{\mcQ} @<{\mu}<< S_{\mcQ} @<{\nu_{z_o}}<<S_{\mcQ, z_o} \\
@V{f'_{\mcQ}}VV                 @VV{f_{\mcQ}}V         @VV{f_{\mcQ, z_o}}V \\
\PP^2@<<{q}< \widehat{\PP^2} @<<{q_{z_o}}< (\widehat{\PP^2})_{z_o},
\end{CD}
\]
where  $q$ is a composition of a finite 
 number of blowing-ups so that
 the branch locus becomes smooth and $q_{z_o}$ is a composition of 
 two blowing-ups. The map
$f_{\mcQ, z_o} : S_{\mcQ} \to (\widehat{\PP^2})_{z_o}$ is the induced double cover which
coincides  with the quotient morphism determined by the involution $[-1]_{\varphi_{\mcQ, z_o}}$ on $S_{\mcQ,z_o}$, which is given by the inversion with respect to the group law 
 on the generic fiber.  
 \item[(iii)] The surface
 $ (\widehat{\PP^2})_{z_o}$ can be blown down to $\Sigma_2$ such that the proper
 transform of $\mcQ$ is mapped to the trisection $\mcT_{S_{\mcQ, z_o}}$  as in 
 Section~\ref{subsec:2-2}. Let $l_{z_o}$ be the tangent line at $z_o$. By our construction,
 $\varphi_{\mcQ, z_o}$ has a singular fiber $F_{z_o}$ containing irreducible componets
 mapped to $l_{z_o}$. The type of $F_{z_o}$ is determined by the residual 
 intersections between $l_{z_o}$ and $\mcQ$. For details on the types of singular
 fibers, see \cite{miranda-persson}.

 \end{enumerate}
 
 Given  $m$-sections $D_i$ and the corresponding sections $s_i:=s_{P_{D_i}}$ ($i = 1, \ldots, r$), we have
 plane curves $\mcC_{D_i}:= f_{\mcQ, z_o}(D_i),\, \mcC_{s_i} (\mbox{ or  $\mcC_{P_{D_i}}$}):= f_{\mcQ, z_o}(s_i)$ $(i = 1, \ldots, r)$. The embedded topology of 
unions of $\mcQ$ and these curves are of
interest. The first and the last authors have previously studied the following cases:

\begin{itemize}

\item In \cite{bannai-tokunaga15}, they considered the case where $\mcQ$ is a union two general conics. 
They  make use of bisections $D_i$ ($i = 1, \ldots r$) such that $\mcC_{D_i}$ $(i = 1, \ldots, r)$ are
conics in order to construct Zariski $N$-tuples for conic arrangements.

\item  In \cite{bannai-tokunaga17}, they consider the case where $\mcQ$ is an irreducible quartic with two nodes. 
They make use of bisections $D_i$ ($i = 1, \ldots, r$) such that $\mcC_{D_i}$ $(i = 1, \ldots, r)$ are
conics in order to give a Zariski $5$-tuple for $\mcQ$ and two conics of this form. 

\item In \cite{bannai-tokunaga20},  they go on to consider trisections of $S_{\mcQ, z_o}$ when $\mcQ$ is a union of 
general four lines. In this case, $E_{S_{\mcQ, z_o}}(\CC(t)) \cong A_1^*\oplus (\ZZ/2\ZZ)^{\oplus 2}$.
In \cite{bannai-tokunaga20}, a generator $P_o$ for the $A_1^*$  part is 
chosen and 
$P_D$ for certain trisections $D$ are explicitly computed in terms of $P_o$ and the torsion part.
As an application, a result given in
\cite{bty20} was refined.

\end{itemize}
See \cite{survey} for definitions of Zariski pairs/triples.

\bigskip

In this article, we continue to study trisections of $S_{\mcQ, z_o}$ which are obtained as follows:

 Let  $E$ be a smooth cubic such that the intersection
 number $i_p(E, \mcQ)$ is even for every $p \in \mcQ\cap E$ and also $z_o\not\in E$. In $S_{\mcQ, z_o}$, $E$ gives rise to either an irreducible curve $\tilde{E}$ which is a 6-section or
 two irreducible components $E^{\pm}$, both of which are trisections of $\varphi_{\mcQ, z_o}$. In the former case, $P_{\tilde{E}}=O$ and  in the latter case, $[-1]P_{E^+} = P_{E^-}$ and 
 $\mcC_{E^\pm}=E$. Note that
 if there exists $p \in \mcQ\cap E$ with $i_p(E, \mcQ)$ being odd, then $E$ gives rise to an an irreducible curve
 in $S_{\mcQ, z_o}$. We apply the methods introduced in \cite{bannai-tokunaga20} in order to compute $P_{E^\pm}$. In fact, we have the following:

 \begin{prop0}\label{prop:trisec-1}{Let $\mcQ$ be a nodal quartic curve, (i.e. a possibly reducible quartic curve whose singularities are at most nodes), which is not a union of four lines. Let $z_o\in\mcQ$ be a general point on a non-linear component of $\mcQ$. Let $E^\pm$ be trisections obtained as above. Then the value of the height pairing of $P_{E^\pm}$ on $S_{\mcQ, z_o}$ is given by
 \[
 \langle P_{E^+}, \, P_{E^+}\rangle =  \langle P_{E^-}, \, P_{E^-}\rangle = \frac 32 - \frac r2,
 \]
 where $r$ is the number of nodes of $\mcQ$ that  $E$ passes through. In particular, $r \le 3$.
 }
 \end{prop0}

We go on to apply this result in the case when  $\mcQ$ consists of components of small degree. We will consider the cases where $\mcQ$ is 
\begin{enumerate}
 \item[(i)]  a union of a smooth conic $C_o$ and two lines $L_1$ and $L_2$ all meeting transversely, or 
\item[(ii)] a union of two smooth conics $C_1$ and $C_2$ meeting transversely,  
\end{enumerate}
and $z_o$ is a smooth point on a conic component. We introduce the following notation for each case:

Case (i):  $L_1\cap L_2 = \{p_0\}$, $L_1\cap C_o = \{p_1, p_2\}$ and $L_2 \cap C_o = \{p_3, p_4\}$.

Case (ii): $C_1\cap C_2 = \{p_1, p_2, p_3, p_4\}$.

%

In these cases, we have
\[
E_{S_{\mcQ, z_o}}(\CC(t)) \cong  \begin{cases}
                                                   (A_1^*)^{\oplus 2}\oplus (\ZZ/2\ZZ)^{\oplus 2}  & \mbox{ for (i)}\\
                                                  (A_1^*)^{\oplus 3}\oplus (\ZZ/2\ZZ)  & \mbox{ for (ii).}
                                                   \end{cases}
 \]
 
As the first and last authors have considered the case where $E_{S_{\mcQ, z_o}}(\CC(t))\cong A_1^\ast\oplus(\ZZ/2\ZZ)^{\oplus2}$,   it seems to  be natural to study the above two cases as a
next step from the viewpoint of $E_{S_{\mcQ, z_o}}(\CC(t))$, and our goal of this article is to give new examples that present new phenomena in the study of the embedded topology of plane curves. 
 
 In what follows, we assume that $\mcQ$ is one of the above (i) and (ii) unless otherwise is stated. 
 We have the following for $\mcC_{P_{E^\pm}}$:
 
 \begin{prop0}\label{prop:trisec-2}{
Under the notation given above, we have:
 \begin{enumerate}
 \item[\rm{(i)}] Suppose that $\mcQ = L_1 + L_2 + C_o$. Then $r = 1, 2$ or $3$. More precisely,  letting
 $\mcP = \{p_0, p_1, p_2, p_3, p_4\}$, we have the following:
     \begin{enumerate}
        \item[\rm{(i-1)}] If $r = 1$, then $E\cap\mcP = \{p_0\}$ and $\mcC_{P_{E^\pm}}$ is a line 
        through $p_0$ and is tangent to $C$.
        
        \item[\rm{(i-2)}] If $r = 2$, then $E\cap\mcP = \{p_i, p_j\},  i \in \{1, 2\}, \, j \in \{3, 4\}$ and 
        $\mcC_{P_{E^\pm}}$ is the line 
        through $p_i$ and $p_j$.
        
         \item[\rm{(i-3)}] If $r = 3$, then $E\cap\mcP = \{p_0, p_1, p_2\}$ or $\{p_0, p_3, p_4\}$  and 
        \[
        \mcC_{P_{E^\pm}} = \left  \{ \begin{array}{cc}
                                               L_1 & \mbox {if $E\cap\mcP = \{p_0, p_1, p_2\}$}, \\
                                               L_2 &  \mbox {if $E\cap\mcP = \{p_0, p_3, p_4\}$}.
                                               \end{array} \right .
       \]
      \end{enumerate}
 \item[\rm{(ii)}]   Suppose that $\mcQ =  C_1 + C_2$. Then $r = 0$ or $2$. More precisely,  letting
 $\mcP = \{p_1, p_2, p_3, p_4\}$, we have the following:
 
    \begin{enumerate}
        \item[\rm{(ii-1)}] If $r = 0$, then  $\mcC_{P_{E^\pm}}$ is a bitangent line to $\mcQ$.
                
        \item[\rm{(ii-2)}] If $r = 2$, then $E\cap\mcP = \{p_i, p_j\}$ and 
        $\mcC_{P_{E^\pm}}$ is the line 
        through $p_i$ and $p_j$.
        
    \end{enumerate}    
             
\end{enumerate}
}
\end{prop0}

  As applications of Proposition~\ref{prop:trisec-2}, we construct Zariski pairs/triples for cubic-conic-line
  arrangements in Section~\ref{subsec:3-3}. In particular,  for the following example based on Proposition~\ref{prop:trisec-2}, (ii-1),
  it may be worthwhile to state it here as a theorem:

 \begin{thm0}\label{thm:trisec-1}{
Let  $\mcQ = C_1 + C_2$ as in  Proposition~\ref{prop:trisec-2} {\rm (ii)}.   There exist two families of plane curves 
  $\{\mcB_{i, s}\}_{s \in \Delta_{\epsilon}} (i = 1, 2)$ over 
$\Delta_{\epsilon} = \{ s \in \CC \mid |s| < \epsilon\}$ for sufficiently small $\epsilon>0$, such that    
  
\begin{enumerate}
  
\item[\rm{(i)}]  For $s \neq 0$, $\mcB_{i, s} = \mcQ + E_{i, s} + L$, $(i =1, 2)$, where $L$ is a fixed bitangent line 
  to $\mcQ$ and  $E_{i,s}$ ($i = 1, 2$) are smooth cubics as in Proposition~\ref{prop:trisec-2} {\rm (ii-1)}
  such that $\mcC_{P_{E^+_{1, s}}} = L$ and $\mcC_{P_{E^+_{2, s}}} \neq L$.
  
\item[\rm{(ii)}] Both of $\mcB_{i, 0}$ $(i = 1, 2)$ are $\mcQ + L_{12} + L_{13}+ L_{23} + L$, where
$L_{ij}$ is the line through $p_i$ and $p_j$.
  
\item[\rm{(iii)}] The pairs $(\mcB_{1, s}, \mcB_{2, s})$  are  Zariski pairs for all $s\not=0$.
\end{enumerate}
In particular there  exists  a family of Zariski pairs which degenerate to the same conic-line arrangement. 	
}
\end{thm0}

\begin{rem}
The existence of such families given in Theorem~\ref{thm:trisec-1} was expected, but the significance of the above Theorem is in giving a simple explicit example. 
\end{rem}


Our new key tool to construct curves as above
is  
^^ the Mumford representation of semi-reduced divisors' on $E_{S_{\mcQ, z_o}}$.
This makes it possible for us to treat multi-sections concretely.
The Mumford representation was first considered in \cite{mumfordII} in order to
describe the Jacobian of  hyperelliptic curves explicitly, based on Jacobi's ideas. It has been exploited
in  the study of hyperelliptic curve cryptography. Our study can be regarded
as a new application of Mumford representations.

 This article consists of $6$ sections.  In Section~1 we first summarize previous results which we need later, and prove Proposition~\ref{prop:trisec-1}. In Section~2, we explain two rational elliptic surfaces which play
 importan roles to prove Proposition~\ref{prop:trisec-2} and Theorem~\ref{thm:trisec-1}. Propositon~\ref{prop:trisec-2} will be proven in Section~3.  We explain the ^^ splitting type' introduced in
 \cite{bannai16} in Section~4. The notion of splitting type plays a key role in order to distinguish the embedded 
 topology of our examples and we prove Theorem~\ref{thm:trisec-1}  (iii)
 assuming the existence of $E_{i, s}$ ($i = 1, 2$).
 In Section~5, we explain the Mumford representation of semi-reduced divisors.  In Section~6, we show the existence of 
 plane curves with the desired properties and construct
 $E_{i, s}$ ($s \neq 0$) in Theorem~\ref{thm:trisec-1} (i) explicitly based on the method given in Section~5.


\section{Elliptic surfaces}\label{sec:elliptic-surfaces}


\subsection{Some terminologies and notation for elliptic surfaces}\label{subsec:1-1}
We here define some notation and  terminologies which are necessary for our later argument. For general
references we refer to \cite{kodaira, miranda-basic, shioda90}. As for bisectios see \cite{bannai-tokunaga15}.
Let $\varphi: S \to C$ be an elliptic surface over a smooth projective curve $C$. Throughout this article,
we always assume that $S$ satisfies the three conditions in the Introduction.
We introduce two subsets $\sing(\varphi)$ and $\Red(\varphi)$ concerning singular fibers as follows, where $F_v = \varphi^{-1}(v)$ for $v \in C$:
\begin{eqnarray*}
\sing(\varphi) &:= & \{ v \in C \mid \mbox{$F_v$ is not a curve of genus $1$.}\}, \\
\Red(\varphi) &:= & \{ v \in \sing(\varphi) \mid \mbox{$F_v$ is reducible.}\}.
\end{eqnarray*}
For $v \in \Red(\varphi)$, $m_v$ denotes the number of irreducible components and we denote the 
irreducible decomposition of $F_v$ by
\[
F_v = \Theta_{v, 0} + \sum_{i=1}^{m_v -1}a_{v, i}\Theta_{v,i},
\]
where $\Theta_{v, 0}$ is the unique component satisfying $\Theta_{v, 0}\cdot O = 1$. We call
$\Theta_{v, 0}$ the identity component. We use Kodaira's symbol (\cite{kodaira}) in order to
describe the types of singular fibers. Irreducible components of singular fibers are labeled
as in \cite[pp. 81-82]{tokunaga14}. 
We denote the set of sections of $\varphi : S \to C$ by $\MW(S)$, which is not empty by our
assumption. The set $\MW(S)$ is endowed with an abelian group structure by considering fiberwise
addition with  $O$ as the zero element. 
Let $E_S$ be the generic fiber of $\varphi : S \to C$. Then $E_S$ can be regarded as a curve of 
genus $1$ over the rational function field, $\CC(C)$, of $C$. Under our setting, $S$ is known as
the Kodaira-N\'eron model of $E_S$ and we can identify 
the set of $\CC(C)$-rational points $E_S(\CC(C))$ of $E_S$  with $\MW(S)$. Here the 
zero of $E_S(\CC(C))$ is the $\CC(C)$-rational point given by the restriction of $O$ to $E_S$. We also
denote it by $O$ for simplicity.

We denote the N\'eron-Severi group of $S$ by $\NS(S)$, which is finitely generated and torsion-free
by \cite[Theorem 1.2]{shioda90}. Following \cite{shioda90}, we define the subgroup $T_{\varphi}$ of $\NS(S)$
generated by $O$, $F$: a general fiber and $\Theta_{v, i}$ $(1\le i \le m_v-1, v \in \Red(\varphi))$.  For a divisor  $D$ on $S$,
the restriction $\fd_D$ of $D$  to $E_S$ is a divisor defined over $\CC(C)$. 
In \cite{shioda90}, Shioda defined a structure on the free part of $E_S(\CC(C))$ called the height pairing,
which we denote by $\langle \bullet, \bullet\rangle$. We use properties of $\langle \bullet, \bullet\rangle$ freely.
For details, see \cite{shioda90}.

 For a reducible singular fiber $F_v$ ($v\in \Red(\varphi)$), $R_v$ denotes the subgroup of $T_{\varphi}$ generated
 by $\Theta_{v, i}$ ($1 \le i \le m_v -1$). We denote its dual by $R^{\vee}_v$. 
The map $\gamma : \MW(S) \to \oplus_{v \in \Red(\varphi)} R^{\vee}_v/R_v$ is the homomorphism given in \cite[p. 70]{miranda-basic}, 
 which describes at which irreducible component a given section meets  $F_v$. This $\gamma$ coincides with the maps considered
 in \cite{bannai-tokunaga20} and  \cite[p. 83]{tokunaga12}.

\subsection{The height pairing  and the homomorphism $\phi_o$}

In this section we review the homomorphism $\phi_o: \Div(S) \to \NS(S)\otimes\QQ$ 
considered in \cite[Section~2]{bannai-tokunaga20} and its relation between 
$\phi_o(D)$ and the height $\langle P_D, P_D\rangle$. Put $\NS(S)_{\QQ}:= \NS(S)\otimes\QQ$.
For $D \in \Div(S)$, we define $\vc(v, D)$ to be
\[          
 \vc(v, D) :=\left [
 \begin{array}{c}
 D\cdot\Theta_{v,1}\\
 \vdots \\
   D\cdot \Theta_{v,m_v-1}
 \end{array} \right ]
 \]
 and put $\FF_v = [\Theta_{v, 1}, \ldots,  \Theta_{v, m_v-1} ]$.
 By \cite[Lemma 5.1]{shioda90}, we have 
 \[
(\ast) \quad D \approx s(D) + (d-1)O + nF +  \sum_{v \in \Red(\varphi)}\FF_v A_v^{-1}(\vc(v, D) - \vc(v, s(D)), 
\]
where $\approx$ denotes algebraic equivalence between divisors, and $d$ and  $n$ 
are integers defined as follows:
\[
d = D\cdot F,  \qquad n = (d-1)\chi({\mathcal O}_S) + O\cdot D - s(D)\cdot O.
\]
Furthermore,  $A_v$ is the intersection matrix $(\Theta_{v,i}\cdot\Theta_{v, j})_{1\le i, j \le m_v-1}$.  
Since  the entries of $A_v^{-1}$ are not necessarily integers, the condition $A_v^{-1}\left (\vc(v, D) - \vc(v, s(D)) \right ) \in \ZZ^{\oplus m_v -1}$ imposes some restriction at which irreducible components of $F_v$, $D$ and $s(D)$
meet.  For example we have the following Lemma:

\begin{lem}\label{lem:I2-fiber}{If $F_v$ is a singular fiber of type $I_2$, $\vc(v, D) - \vc(v, s(D))$ is 
even (Note that $\vc(v, D)$ becomes an integer in this case).
}
\end{lem}

Let $\phi_o : \Div(S) \to \NS(S)_{\QQ}$ be the homomorphism  introduced in \cite{bannai-tokunaga20} whose explicit form is
\[
 (\ast\ast) \quad  \phi_o(D) = D - dO - (d\chi({\mathcal O}_S) + O\cdot D)F - \sum_{v \in \Red(\varphi)}\FF_v A_v^{-1}\vc(v, D).
 \]
Let $\phi: E_S(\CC(C)) \to \NS(S)_{\QQ}$ be the  homomorphism given in \cite[Lemma 8.2]{shioda90}. Then we have

\begin{lem}\label{lem:fundamental}{ {\rm \cite[Lemma2.1]{bannai-tokunaga20}}
 \begin{enumerate}
 
\item[{\rm (i)}] For $P \in E_S(\CC(C))$ and its corresponding section $s_P$, we have $\phi(P) = \phi_o(s_P)$.
 
\item[{\rm (ii)}] For $D \in \Div(S)$ and its corresponding point $P_D \in E_S(\CC(C))$,  we have
$\phi(P_D) = \phi_o(D)$.
 
\end{enumerate}
}
\end{lem}

\begin{cor}\label{cor:fundamental}{For divisors $D_1, D_2$, we have
\[
\langle P_{D_1}, P_{D_2} \rangle = - \phi_o(D_1)\cdot\phi_o(D_2).
\]
}
\end{cor}

\proof As $\langle P_{D_1}, P_{D_2} \rangle = - \phi(P_{D_1})\cdot\phi(P_{D_2})$, our statement is 
immediate from Lemma~\ref{lem:fundamental}. \endproof

By Corollary~\ref{cor:fundamental}, $\langle P_{D_1}, P_{D_2} \rangle$ can be computed by the geometric
data of $D_1$ and $D_2$. Also, by calculating $\phi_o(D_1)\cdot\phi_o(D_2)$ using $(\ast\ast)$, we have
\[
 (\ast\ast\ast) \quad\phi_o(D_1) \cdot\phi_o(D_2)=D_1\cdot D_2 - d_2D_1\cdot O - d_1D_2\cdot O -d_1d_2\chi({\mathcal O}_S)-\sum_{v\in\Red(\varphi)} {}^t\vc(v, D_1)A_v^{-1}\vc(v, D_2),
\]
where $d_i = D_i\cdot F$ $(i = 1, 2)$.
In particular, if $D_1 = D_2$, we have
\[
\phi_o(D_1) \cdot\phi_o(D_1)=D_1^2-2d_1D_1\cdot O-d_1^2\chi({\mathcal O}_S)-\sum_{v\in\Red(\varphi)} {}^t\vc(v, D_1)A_v^{-1}\vc(v, D_1).
\]

Furthermore we have the following lemma which we will use later.
\begin{lem}\label{lem:int-trisec}
If the reducible singular fibers of $S$ are  all of type $\I_2$ or $\III$ and both $D_1$ and $D_2$ are  integral 
$d_1$- and $d_2$-sections respectively, then we have
\[
\langle P_{D_1}, P_{D_2} \rangle =-\phi_o(D_1)\cdot\phi_o(D_2)=-\left(D_1\cdot D_2 -d_1d_2\chi(\mathcal{O}_S)+\frac{1}{2}\sum_{v\in\Red(\varphi)}{}^t\vc(v,D_1)\vc(v,D_2)\right).
\] 
\end{lem}
\proof Since the reducible singular fibers are all of type $\I_2$, we have $A_v^{-1}=\left[-\frac{1}{2}\right]$ for all $v\in \Red(\varphi)$. Also, since $D_i$ (i = 1, 2) are  integral trisections, we have $D_i\cdot O=0$  and $d_i=D_i\cdot F=3$ $(i = 1, 2)$. Then our result is obtained directly from $(\ast\ast\ast)$.
\endproof

 \subsection{The proof of Proposition \ref{prop:trisec-1} }\label{subsec:1-3}

In this subsection we prove Proposition \ref{prop:trisec-1}. Let $\mcQ$ be a nodal quartic curve which is not a union of four lines and let $z_o\in \mcQ$ be a general point on a non-linear component of $\mcQ$. By these assumptions, $S_{\mcQ, z_o}$ is a rational elliptic surface whose reducible singular fibers are all of type $\I_2$ or $\III$. Note that one reducible singular fiber $F_\infty$ corresponds to the tangent line $L_\infty$ of $\mcQ$ at $z_o$ and the other singular fibers correspond to the lines passing through $z_o$ and a node of $\mcQ$.

Let $E$ be a smooth cubic curve such that the intersection
 number $i_p(E, \mcQ)$ is even for every $p \in \mcQ\cap E$ and $z_o\not\in E$. We assume that $E$ give rises to two trisections $E^\pm$ in $S_{\mcQ, z_o}$. Note that since $z_o\not\in E$, $E^\pm$ become integral trisections. Let $r$ be the number of the nodes of $\mcQ$ that $E$ passes through. Let $\widehat{E}\subset (\widehat{\PP})_{z_o}$ be the strict transform of $E$ under  the compositions of blow-ups $q_{z_o}\circ q$. Then we have
$\widehat{E}\cdot \widehat{E}=9-r$ since $E$ is blown-up once at each of the $r$ nodes.  Furthermore, since $f_{\mcQ,z_o}^\ast(\widehat{E})=E^++E^-$, we have
\[
(E^++E^-)^2=2\widehat{E}\cdot \widehat{E}=2(9-r)
\]  
 and since $(E^+)^2=(E^-)^2$, we obtain
 \[
 (E^+)^2+E^+\cdot E^-=9-r.
 \]
 On the other hand, since $E^+\cdot E^-=\frac{1}{2}\widehat{\mcQ}\cdot \widehat{E}$, where $\widehat\mcQ$ is the strict transform of $\mcQ$ in  $S_{\mcQ, z_o}$, we have $E^+\cdot E^-=6-r$. Hence we have $(E^\pm)^2=3$. 
 
 Next, we consider ${}^t\vc(v,E^\pm)$. Let $F_\infty=\Theta_{\infty,0}+\Theta_{\infty,1}$, where $\Theta_{\infty,0}$ is the pre-image under $f_{\mcQ, z_o}$ of the exceptional divisor of the last blow-up in  $q_{z_o}\circ q$,  and $\Theta_{\infty,1}$ is the pre-image of the strict transform of the tangent line $L_\infty$. Then since $E\cdot L_\infty=3$, we have $\vc(\infty,E^\pm)=[3]$ and ${}^t\vc(\infty,E^\pm)\vc(\infty,E^\pm)=9$. Next, let $F_v$ be a reducible fiber corresponding to a node $v$ of $\mcQ$. Then $F_v=\Theta_{v,0}+\Theta_{v,1}$, where $\Theta_{v,0}$ is the pre-image under $f_{\mcQ, z_o}$of the strict transform $\widehat{L}_v$ of the line $L_v$ through $z_o$ and $v$, and $\Theta_{v,1}$ is the pre-image $\widehat{E}_v$ of the  strict transform of the exceptional divisor $E_v$ of the blow-up at $v$. Then since 
 \[
 \widehat{E}\cdot \widehat{E}_v=\begin{cases} 1 & ( v\in E) \\
 0 & (v\not\in E) 
 \end{cases}
 \] 
 we have
 \[
 \vc(v,E^\pm)=\begin{cases}
 1 & (v\in E) \\
 0 & (v\not\in E).
 \end{cases}
 \]
 Now since all of the reducible singular fibers of $S_{\mcQ, z_o}$ are of type $\I_2$ and $E^\pm$ are integral trisections, we can apply Lemma \ref{lem:int-trisec} and obtain 
\begin{align*}
\langle P_{E^+}, P_{E^+} \rangle =\langle P_{E^-}, P_{E^-}\rangle &=-\left((E^+)^2-9\chi(\mathcal{O}_{S_{\mcQ, z_o}})+\frac{1}{2}\sum_{v\in\Red(\varphi)} {}^t\vc(c,E^+)\vc(c,E^+)\right)\\
&=-\left(3-9+\frac{1}{2}\left(9+r\right)\right)\\
&=\frac{3}{2}-\frac{1}{2}r
\end{align*}
which proves Proposition \ref{prop:trisec-1}.

\section{Two rational elliptic surfaces for the cases (i) and (ii)}\label{sec:two-elliptic surfaces}




 In this section, we give  rather detailed descriptions of rational elliptic surfaces for the cases (i) and (ii) in the
 Introduction. Both of them were considered in \cite{tokunaga14} and explicit models for both cases were given.  As for the translation by a $2$-torsion,
see \cite[Lecture VII]{miranda-basic} or \cite{tokunaga12}.

\subsection{The case (i): A rational elliptic surface attached to a conic and two lines}\label{subsec:case-1}

Let $\mcQ = C_o + L_1 + L_2$ be as in the Introduction and we keep our notation there.
Take a smooth point of $C_o$ as $z_o$  and let 
$ S_{\mcQ, z_o}$ be the rational elliptic
surface associate to $\mcQ$ and $z_o$ and we denote its elliptic fibration by
$\varphi_{\mcQ, z_o} : S_{\mcQ, z_o} \to \PP^1$. As we have seen in \cite{tokunaga14},
$S_{\mcQ, z_o}$ satisfies the following:
\begin{enumerate}
\item[(i)] The elliptic fibration $\varphi_{\mcQ, z_o}$ has $6$ singular fibers of type $\I_2$ or $\III$ and no other
singular fiber. They arise from the tangent line at $z_o$ and lines connecting $z_o$ and $p_i$ $(i = 0, 1, \ldots, 4)$. We denote them by $F_{\infty}$ and $F_i$ $(i = 0, 1, \ldots, 4)$, respectively and their irreducible decompositions by $F_{\bullet} = \Theta_{\bullet, 0} + \Theta_{\bullet, 1}$
$(\bullet = \infty, 0, 1,\ldots, 4)$.

\item[(ii)] The group $E_{S_{\mcQ, z_o}}(\CC(t))$ is isomorphic to  $(A_1^*)^{\oplus 2}\oplus (\ZZ/2\ZZ)^{\oplus 2}$. 

\item[(iii)] Let $L_{ij} (i<j)$ be the line connecting $p_i$ and $p_j$. Let $L_1 = L_{12}$ and 
$L_2=L_{34}$. The $6$ lines $L_{ij}$ give rise to sections of $\varphi_{\mcQ, z_o}$ which can be identified with points in $E_{S_{\mcQ, z_o}}(\CC(t))$.
More precisely, both of $L_1$ and $L_2$ give  $2$-torsion points, and as for the other
four lines $L_{13}, L_{14}, L_{23}$ and $L_{24}$,  each of them give rise to two points, resulting in eight points
$P_{ij}, [-1]P_{ij}$ $(i, j) = (1, 3), (1,4)$, $(2, 3), (2, 4)$. We may assume
\[
(A_1^*)^{\oplus 2} \cong \ZZ P_{13}\oplus \ZZ P_{14}
\]
and that $P_{23}, P_{24}$ are obtained by  translation by suitable $2$-torsion points.

\item[(iv)] For each $P_{ij}$, $\mcC_{[2]P_{ij}}$ is a conic inscribed by  $\mcQ$.
\end{enumerate}


\subsection{The case (ii): A rational elliptic surface attached to two coincs}\label{subsec:case-2}

\subsubsection{General properties for case (ii)}

Let $\mcQ = C_1 + C_2$ be as in the Introduction and  we keep the notation in the Introduction.
  Take
a smooth point $z_o$ of $\mcQ$ such that $l_{z_o}$ is not a bitangnet line to $\mcQ$. The rational elliptic surface $\varphi_{\mcQ, z_o}: S_{\mcQ, z_o} \to \PP^1$ associate to $\mcQ$ and $z_o$ satifies the following properties:

\begin{enumerate}
\item[(i)] The elliptic fibration $\varphi_{\mcQ, z_o}$ has $5$ singular fibers of type $\I_2$ or $\III$ and no other
reducible singular fiber. They arise from the tangent line at $z_o$ and lines connecting $z_o$ and $p_i$ $(i = 1, \ldots, 4)$.
We denote them by $F_{\infty}$ and $F_i$ $(i = 1, \ldots, 4)$, respectively and their irreducible decompositions by $F_{\bullet} = \Theta_{\bullet, 0} + \Theta_{\bullet, 1}$
$(\bullet = \infty, 1,\ldots, 4)$.

%

\item[(ii)] The group $E_{S_{\mcQ, z_o}}(\CC(t))$ is isomorphic to  $(A_1^*)^{\oplus 3}\oplus \ZZ/2\ZZ$. 
We denote
the unique $2$-torsion point by $T$, which arises from the conic containing $z_o$. 

\item[(iii)] Let $L_{ij}$ $(i < j)$ denote a line connecting $p_i$ and $p_j$. The $6$ lines
$L_{ij}$ ($1\le i < j \le 4$) give rise to $12$ rational points $P_{ij}$ and $[-1]P_{ij}$
($1 \le i < j \le 4$) of $E_{S_{\mcQ,z_o}}(\CC(t))$, each of which has the height
$\langle P_{ij}, P_{ij} \rangle = 1/2$. We may assume
\[
(A_1^*)^{\oplus 3} \cong \ZZ P_{12} \oplus \ZZ P_{13} \oplus \ZZ P_{23}
\]
and  $P_{ij} \dot {+} T = P_{kl}$, where $\{i, j, k, l\} =\{ 1, 2, 3, 4\}$.

\end{enumerate}

\subsubsection{Rational points arising from lines in case (ii)}

We consider the structure of  $S_{\mcQ, z_o}$ in more detail for later use. 
 
 \begin{lem}\label{lem:line-sec1}{The curve $\mathcal{C}_P$ is a line if and only if $s_P\cdot\Theta_{\infty, 1} = 1$ and $s_P\cdot O = 0$.
 In particular, if $\mathcal{C}_P$ is a line, then $z_o \not\in \mathcal{C}_P$.
 }
 \end{lem}
 
 \proof Choose a general fiber $F$ of $\varphi_{\mcQ, z_o}$. By our construction, $q\circ q_{z_o}\circ f_{\mcQ, z_o}(F)$ is
 a line through $z_o$. We may assume $F\cap O \cap s_P  = \emptyset$.  If $\mathcal{C}_P$ is a line, $\mathcal{C}_P$ and 
 $q\circ q_{z_o}\circ f_{\mcQ, z_o}(F)$ meet one point, which is the image $z_1$ of $s_P\cap F$. By our choice of $F$,
 $z_1 \neq z_o$ and $z_o \not\in \mathcal{C}_P$.  As both $O$ and $\Theta_{\infty, 0}$ are mapped to $z_o$, $s_P\cdot \Theta_{\infty, 1}
 = 1$ and $s_P\cdot O = 0$. Conversely, if $s_P\cdot\Theta_{\infty,1} =1$ and $s_P\cdot O = 0$, then 
 $\mcC_P\cap q\circ q_{z_o}\circ f_{\mcQ, z_o}(F) = q\circ q_{z_o}\circ f_{\mcQ, z_o}(s_P\cap F)$.
 As $\sharp(s_P\cap F) = 1$,  $\mathcal{C}_P$ is a line.
 \endproof
 
 \begin{lem}\label{lem:line-sec2}{Let $P \in E_{S_{\mcQ, z_o}}(\CC(t))$ be a point such that $\mathcal{C}_P$ is a line. 
 \begin{enumerate}
 \item[\rm (i)]  For $p \in \mathcal{C}_P\cap \mcQ$, the intersection multiplicity $i_p(\mcC_P, \mcQ)$ at $p$ is even.
 \item[\rm (ii)] The height pairing of $P$ is given by
 \[
 \langle P, P \rangle = \frac 32 - \frac 12\sharp(\mcC_P\cap C_1\cap C_2).
 \]
 
 \end{enumerate}
 Conversely, if $l$ is a line such that {\rm (a)} $z_o \not\in l$ and  {\rm (b)}  $i_p(l, \mcQ)$ is even for $\forall p \in l\cap \mcQ$, then
 $l$ gives rise to two points $P_l$ and $[-1]P_l$ in $E_{S_{\mcQ, z_o}}(\CC(t))$ such that 
 \[
  \langle P_l, P_l \rangle = \frac 32 - \frac 12\sharp(l \cap C_1\cap C_2).
 \]
 }
 \end{lem}
 
 \proof (i)  If there exists a point $p \in \mcC_P\cap Q$ such that
 $i_p(\mcC_P, \mcQ)$ is odd, $f'_{\mcQ, z_o}$ gives a ramified double cover of $\mcC_P$, this is impossible as the image of $s_{[-1]P}$ is also $\mcC_P$. 
 
 (ii) Since $s_P\cdot O = 0$,  $s_P\cdot\Theta_{z_o, 1} = 1$ and  $s_P\cdot \Theta_{i, 1} = 1$ if and only if $p_i \in \mcC_P$. Hence
 by \cite[Theorem 8.6, (8.12)]{shioda90}, our statement follows. The remaining statements follows easily as $\mcC_{P_l} = l$.

\endproof
 
 \subsubsection{Bitangents to $\mcQ$ and triangles $\triangle_{ijk} = L_{ij} + L_{ik} + L_{jk}$ in case (ii)}
 \label{subsubsec:bitangent}
  
 Let $s_T$ be the section corresponding to $T$.   Assume that $z_o \in C_1$. Then we infer that
 $C_1$ gives rise to a $2$-torsion point $P_{C_1}$, i.e., $P_{C_1} = T$. Hence we have
 $s_T\cdot\Theta_{i, 1} = 1$ and $s_T\cdot\Theta_{\infty, 0} = 1$. Thus we see $P_{ij}\dot{+}T = P_{kl}$ by the property of the
 map $\gamma$ (see \cite[p. 70]{miranda-basic}. Also \cite[Theorem 9.1]{kodaira} for explicit description).
 
  We now consider rational points in $E_{S_{\mcQ, z_o}}(\CC(t))$ given by a triangle 
  $\triangle_{123} := L_{12} + L_{13} + L_{23}$.
  Since each $L_{ij}$ produces $P_{ij}$ and $[-1]P_{ij}$  ($[\pm 1]P_{ij}$ for short), 
   we have $8$ rational points in $E_{S_{\mcQ, z_o}}(\CC(t))$ from $\triangle_{123}$:
\begin{eqnarray*}
 Q_1 &:= &  P_{12}\dot{+}P_{13}\dot{+}P_{23},  \\
 Q_2 &:= &  [-1]P_{12}\dot{+}P_{13}\dot{+}P_{23}, \\ 
  Q_3 &:= & P_{12}\dot{+}[-1]P_{13}\dot{+}P_{23},  \\
  Q_4 &:= & P_{12}\dot{+}P_{13}\dot{+}[-1]P_{23}, 
 \end{eqnarray*}
 and $[-1]Q_j$ $(j = 1, 2, 3, 4)$. As $P_{12}, P_{13}, P_{23}$ is a basis of the 
 free part of $E_{S_{\mcQ, z_o}}(\CC(t))$, these $8$ points are distinct. For
 these points, we have
 \begin{itemize}
 \item $\langle Q_i, Q_i \rangle = \langle [-1]Q_i, [-1]Q_i\rangle =  3/2$ and
 \item $\mcC_{Q_i} = \mcC_{[-1]Q_i}$.
 \end{itemize}
 By the property of $\gamma$, we infer that $s_{Q_j}\cdot \Theta_{i,1} = s_{[-1]Q_j}\cdot
 \Theta_{i,1} = 0$ $(j = 1, 2, 3, 4)$ and 
 $s_{Q_j}\cdot \Theta_{\infty,1} = s_{[-1]Q_j}\cdot
 \Theta_{\infty,1} = 1$ $(j = 1, 2, 3, 4)$. Hence by \cite[Theorem 8.6, (8.12)]{shioda90},
 $s_{Q_j}\cdot O = s_{[-1]Q_j}\cdot O = 0$ $(j = 1, 2, 3, 4)$. Thus by 
 Lemmas~\ref{lem:line-sec1} and \ref{lem:line-sec2}, $\mcC_{Q_j}$ ($j= 1, 2, 3, 4$) are the
 four bitangent lines of $\mcQ$. For other triangles $\triangle_{ijk}$ $(1\le i< j < k \le 4)$,
 as $P_{ij}\dot{+}T = P_{kl}$ for $\{i, j, k,l\} = \{1, 2, 3, 4\}$, we have the same rational
 points $[\pm 1]Q_j$ ($j = 1, 2, 3, 4$) and the same $4$ bitangent lines.

\section{The proof of Proposition~\ref{prop:trisec-2}}

We prove the statements for $\mcC_{E^+}$ only, since $P_{E^-} = [-1]P_{E^+}$ and $\mcC_{P_{E^-}}
= \mcC_{P_{E^+}}$.
As $\vc(i, E^+) = \vc(i, s_{P_{E^+}})$ $(i = 0, 1,\ldots, 4)$ and $\vc(\infty, s_{P_{E^+}}) = [1]$, we have:
\[
\langle P_{E^+}, P_{E^+} \rangle 
  =  \frac {3}{2} - \frac {1}{2} r 
\]
by Proposition~\ref{prop:trisec-1}. On the other hand, by the explicit formula of the height pairing, we have
 \[ \langle P_{E^+}, P_{E^+} \rangle =  \frac {3}{2} + s_{P_{E^+}}\cdot O - \frac {1}{2} r,\]
hence $s_{P_{E^+}}\cdot O = 0$. Now, by Lemma~\ref{lem:line-sec1}, $\mcC_{P_{E^+}}$ is a line not passing
through $z_o$.

\subsection{Proof for the case (i) $\mcQ = C_o + L_1 + L_2$}\label{subsec:prop3-case1}

We first show that $r \neq 0$. If $r = 0$, there exists $p \in L_1\setminus \mcP$ such that
$i_p(E, L_1) = i_p(E, \mcQ) = 1$. This is impossible as $i_p(E, \mcQ)$ is even for all $p \in E\cap \mcQ$.


\begin{itemize}
\item[(i-1)]  Suppose that $r = 1$. If $\mcC_{P_{E^+}}\cap \mcP = \{p_1\}$, then $E$ meets $L_2$ at
 $L_2\setminus \mcP$. This means that $E\cap L_2$ contains a point $p$ with $i_p(E, \mcQ) = 1$, 
 which is impossible as $E$ gives rise to $E^{\pm}$. Similarly $p_i \not\in\mcC_{P_{E^+}}\cap \mcP$.
 Hence $\mcC_{P_{E^+}}\cap \mcP = \{p_0\}$. Now since $\vc(\infty, s_{P_{E^+}}) = \vc(0, s_{P_{E^+}}) = [1]$,
 $\vc(i, s_{P_{E^+}}) = [0]\,\,(i = 1, \ldots, 4)$, we infer that $\mcC_{P_{E^+}}$ is a line through $p_0$ and 
 tangent to $C_o$.
 
\item[(i-2)] Suppose that $r=2$. If $p_0 \in E\cap \mcQ$, then by a similar argument to 
the case (i-1), we infer that
 $E$ has a transversal intersection point with either $L_1$ or $L_2$, which is not a node of $\mcQ$.
 Hence $p_0 \not\in E \cap \mcQ$.  If $p_0 \in E\cap \mcQ =\{p_1, p_2\}$ or $\{p_3, p_4\}$, we infer
 that $E$ must passes through $p_0$. Hence these cases do not occur and $\mcC_{P_{E^+}}$ is
 a line described in the statement.
 
\item[(i-3)] As we see in (i-2), If $\{p_1, p_2\}$ or $\{p_3, p_4\} \subset E \cap \mcQ$, then 
 $p_0 \in  E \cap \mcQ$.  Hence $\mcC_{P_{E^+}}$ is a line described in the statement.
\end{itemize}
\subsection{Proof for the case (ii) $\mcQ = C_1 + C_2$}\label{subsec:prop3-case2}

If $r =1$ or $3$, we infer that $E$ intersects $\mcQ$ at a smooth point $p$ with
$i_p(E, \mcQ)$ odd.
 This
is impossible as $i_p(E, Q)$ is even for any $p \in E\cap \mcQ$. Hence $r = 0$ or $2$. 
Since $\vc(i, E^+) = \vc(i, s_{P_{E^+}})$, we infer that $\mcC_{P_{E^+}}$ is a line described in 
the statement.

\section{Splitting types and the topology of curve arrangements}

In this section, we recall the basics  on \textit{splitting types} introduced in \cite{bannai16} which can be used to distinguish the embedded topology of certain plane curves. Furthermore, we calculate the splitting types for certain plane curves derived from sections and trisections of rational elliptic surfaces considered in 
the previous section.

\subsection{The definition of splitting types}

Let $B\subset \PP^2$ be a reduced curve of even degree and let $\pi_B: B\rightarrow \PP^2$ be  the double cover of $\PP^2$ branched along $B$. An irreducible  plane curve $C\subset \PP^2$ is called a {\it contact curve} of $B$ if the local intersection multiplicities $i_{p}(B, C)$ are even for $\forall p\in B\cap C$.  If $C$  is a contact curve, the pre-image $\pi_B^\ast(C)$ of $C$  under $\pi_B$ may be reducible, in which case we call $C$ a {\it splitting curve} with respect to $B$.

\begin{ex}{\rm
Each of the 6 lines $L_{ij}$ in  Section~\ref{subsec:case-1} are splitting curves with respect to $\mcQ$. 
}
\end{ex}


Now given a pair of splitting curves $C_1, C_2$ with respect to $B$, we define the splitting type as follows. Let $\pi_B^\ast(C_i)=C_i^++C_i^-$ be the irreducible decomposition of $\pi_B^\ast(C_i)$ ($i=1,2$). 

\begin{defin}
Let $m_1\leq m_2$ be non-negative integers. Suppose that $C_{1}\cap C_{2}\cap B=\emptyset$. The splitting type of the triple $(C_1, C_2; B)$ is defined to be $(m_1, m_2)$ if for a suitable choice of labels of $C_1^\pm, C_2^\pm$ the equalities $C_1^+\cdot C_2^+=m_1$ and $C_1^+\cdot C_2^-=m_2$ hold.
\end{defin}

The following proposition allows us to distinguish the embedded topology of cures of the form $B+C_1+C_2$ by studying their splitting types.

\begin{prop}\label{prop:splitting-type}[\cite{bannai16}, Proposition 2.5]
Let $B_1, B_2$ be  plane curves of  even degree and let $C_{j1}, C_{j2}$ be splitting curves with respect to $B_j$ $(j=1,2$). Suppose that $C_{j1}\cap C_{j2}\cap B_j=\emptyset$, $C_{j1}$ and $C_{j2}$ intersect  transversally and that $(C_{11}, C_{12}; B_1)$ and $(C_{21}, C_{22}; B_2)$ have distinct splitting types. Then a homeomorphism $h:\PP^2\rightarrow\PP^2$ such that $h(B_1)=B_2$ and $\{h(C_{11}), h(C_{12})\}=\{C_{21}, C_{22}\}$ does not exist.
\end{prop}

\subsection{Computations of splitting types for some curve arrangements related to Proposition~\ref{prop:trisec-1}}

We keep the previous notation given so far.

\subsubsection{The case (i)}

 Let $\mcQ = C_o + L_1 + L_2$ be a quartic as in Proposition~\ref{prop:trisec-1} and let 
 $\varphi_{\mcQ, z_o} : S_{\mcQ, z_o} \to \PP^1$ be the rational elliptic surface in Section~\ref{subsec:case-1}.
 Let $E_1$ and $E_2$ be smooth cubics as in Proposition~\ref{prop:trisec-2} (i-2) such that
 $E_1\cap\mcP =\{p_1, p_3\}$ and $E_2\cap \mcP = \{ p_1, p_4\}$. Let $E_1^{\pm}$ and $E_2^{\pm}$ be
 trisections arising from $E_1$ and $E_2$, respectively. By Proposition~\ref{prop:trisec-2} (i-2) , we may
 assume that $P_{E_1^+} = P_{13}$ and $P_{E_2^+} = P_{14}$. Note that $E_i^{\pm}$ $(i = 1,2)$ are
 all integral trisections and 
 \[
 \vc(\infty, E_1^+) = [3], \, \vc(i, E_1^+) = [0],\, (i = 0, 2, 4), \, \vc(i, E_1^+) = [1] \, (i = 1, 3)
 \]
 and 
  \[
 \vc(\infty, E_2^+) = [3], \, \vc(i, E_2^+) = [0],\, (i = 0, 2, 3), \, \vc(i, E_2^+) = [1] \, (i = 1, 4).
 \]
 Let $\mcC_{[2]P_{13}}$ be the conic as in Section~\ref{subsec:case-1}. We compute the splitting types of
 $(E_1, \mcC_{[2]P_{13}}; \mcQ)$ and $(E_2, \mcC_{[2]P_{13}}; \mcQ)$. Put 
 $f_{\mcQ, z_o}^*\mcC_{[2]P_{13}} = C^+ + C^-$.
 By Lemma~\ref{lem:int-trisec}, we have
 \[
 \langle P_{E_1^+}, P_{C^{\pm}} \rangle = - (E_1^+\cdot C^{\pm} - 3).
 \]
 As $\langle P_{E_1^+}, P_{C^{+}} \rangle = 1$ and $\langle P_{E_1^+}, P_{C^{-}} \rangle = -1$, the splitting
 type of $(E_1, \mcC_{[2]P_{13}}; \mcQ)$ is $(2, 4)$. 
 A similar computation shows that the splitting
 type of $(E_2, \mcC_{[2]P_{13}}; \mcQ)$ is $(3, 3)$. 
 
 \subsubsection{The case (ii)}
 
 Let $\mcQ = C_1 + C_2$ be a quartic as in Proposition~\ref{prop:trisec-1} and let 
 $\varphi_{\mcQ, z_o} : S_{\mcQ, z_o} \to \PP^1$ be the rational elliptic surface in Section~\ref{subsec:case-2}.
 Let $Q_1$ and $Q_2$ be the rational points as in Section~\ref{subsubsec:bitangent}. Note that $\mcC_{Q_i}$
 $(i = 1, 2)$ are bitangents to $\mcQ$. Choose two smooth cubics $E_1$ and $E_2$ such that
 the trisections $E_i^{\pm}$ arising from $E_i$ satisfies $P_{E_i^+} = Q_i$ for each $i$.
 Note that $E_i^{\pm}$ $(i = 1,2)$ are
 all integral trisections and 
 \[
 \vc(\infty, E_1^+) = [3], \, \vc(i, E_1^+) = [0],\, (\forall i),
 \]
 and 
  \[
 \vc(\infty, E_2^+) = [3], \, \vc(i, E_2^+) = [0],\, (\forall i).
 \]
  We now compute
 the splitting types of $(E_1, \mcC_{Q_1}; \mcQ)$ and $(E_2, \mcC_{Q_1}; \mcQ)$.
 Put $f^*_{\mcQ, z_o}(\mcC_{Q_1}) = s_{Q_1} + s_{[-1]Q_1}$.
 Since $Q_1 = P_{12}\dot{+}P_{13}\dot{+}P_{23}$ and $Q_2 = [-1]P_{12}\dot{+}P_{13}\dot{+}P_{23}$, we have
 \[
 \langle P_{E_1^+}, [\pm1] Q_1\rangle = \langle Q_1, [\pm1]Q_1\rangle = \pm \frac 32,
 \]
 and 
 \[
 \langle P_{E_2^+}, [\pm1]Q_1\rangle = \langle Q_2, [\pm1]Q_1\rangle = \pm \frac 12.
 \]
By Lemma~\ref{lem:int-trisec}, we have
\[
 \langle P_{E_i^+}, [\pm 1]Q_1 \rangle = - \left (E_i^+\cdot s_{[\pm1]Q_1} - \frac 32 \right), \, (i = 1, 2).
 \]
Therefore the splitting type of $(E_1, \mcC_{Q_1}; \mcQ)$ is $(0, 3)$, while that of 
 $(E_2, \mcC_{Q_1}; \mcQ)$ is $(1, 2)$.

\subsection{Some Zariski pairs and the proof of  Theorem~\ref{thm:trisec-1} (iii)}\label{subsec:3-3}

\begin{itemize}
\item {\sl A Zariski pair from the case (i).} Let $\mcQ$, $E_i\,(i = 1, 2)$ and $\mcC_{[2]P_{13}}$ be those in 
the case (i) as above. Put $\mcB_i = \mcQ + \mcC_{[2]P_{13}} + E_i$ $(i = 1, 2)$. If 
$\mcB_1$ and $\mcB_2$ have the same combinatorics, then by Proposition~\ref{prop:splitting-type},
$(\mcB_1, \mcB_2)$ is a Zariski pair. We give an explicit example in Section~\ref{subsec:6-1}.

\item {\sl A Zariski pair from the case (ii) and the proof of  Theorem~\ref{thm:trisec-1} (iii).}
Assume that $E_1$ and $E_2$ are as in the case (ii) as above.
Put $\mcB_i = \mcQ + E_i + \mcC_{Q_1}$ $(i = 1, 2)$. If $\mcB_1$ and $\mcB_2$
have the same combinatorics, then $(\mcB_1, \mcB_2)$ is a Zariski pair.
As for the existence of families of cubics $E_{i, s}$, $s \in \Delta_{\epsilon}$ such that (i) for $s \neq 0$, 
$E_{1, s}$ and $E_{2, s}$ satisfy the condition on the cubics $E_1$ and $E_2$ as above, 
and (ii) $E_{1, 0} = E_{2, 0} = L_{12} + L_{13} + L_{23}$,
it is shown in Section~\ref{subsec:6-2}.
\end{itemize}

\begin{rem}\label{rem:triple}{
So far, we have only considered smooth cubics $E$ such that it gives rise to two components. Another case, i.e.,
$E$ gives an irreducible component $\tilde E$, actually happens for the cases (i) and (ii). Let $E_3$ be such a cubic.
By a similar argument to that of \cite{bannai-tokunaga20-1}, we see that $(\mcQ + E_j, \mcQ + E_3)$ $(j = 1, 2)$
are Zariski pairs if $\mcQ+ E_j$ $(j= 1, 2)$ and $\mcQ + E_3$ have the same combinatorics. Explicit construction $E_3$
can be found in \cite{takatoku20-1}. Thus we have a Zariski triple for each case.
}
\end{rem}

\begin{rem}\label{rem:dihedral}{We also note that the above Zariski pairs can be distinguished by studying the existence/non-existence of certain dihedral covers branched along $\mcB_i$ ($i=1,2$), as considered in \cite{tokunaga14}.
}
\end{rem}

\section{Mumford Representations}\label{sec:mumford-representation}


\subsection{Semi-reduced divisors on hyperelliptic curves and their representations}

For terminologies for curves defined over a field $K$, we refer to \cite{silverman}.
Let $C$ be a hyperelliptic curve of genus $g$  defined over $K$ (${\rm char}\,K\not=2$) given by an
affine equation
\[
y^2 = f(x) = x^{2g+1} + \ldots + c_{2g+1} \quad c_i \in K\quad  (i =1, \ldots, 2g+1).
\]
We denote the point at infinity by $O$. By considering $\pi: (x, y) \mapsto x$,
$C$ can be considered as a double cover of $\PP^1$. We denote the covering morphism by
$\pi : C \to \PP^1$ and the hyperelliptic involution by $\iota : (x, y) \mapsto (x, -y)$. The branch points
of $\pi$ are the zeros of $f(x) = 0$ and $x = \infty$ over which we have
$O$. Let $\Div^0(C)$  and $\pic^0(C)$ be the group of divisors  and the divisor class group
of degree $0$ on $C$, respectively. In \cite{mumfordII}, Mumford gave a description of an element of 
$\pic^0(C)$ by two polynomials in $\Kbar[x]$, based on Jacobi's idea. It has been used  in the study of hyperelliptic curve
cryptography in order to compute the addition on $\pic^0(C)$. 
We here explain the Mumford representation briefly.

For a divisor $\fd = \sum_{P\in C}m_P P$ on $C$, we put $\Supp(\fd):= \{P\in C \mid m_P \neq 0\}$. A divisor $\fd$ is
said to be an affine divisor if $O \not\in \Supp(\fd)$.
For an effective affine divisor $\fd$, we can rewrite $\fd$ in the form
\[
\fd = \fd_{\sr} + \pi^*\fd_o, \quad \fd_{\sr} = \sum_{P \in C} n_P P \quad (n_P\geq 0), 
\]
where $\fd_o$ is a divisor on $\PP^1$ and 
\begin{itemize}

\item if   $P \in \Supp(\fd)$ and $P\neq \iota(P)$, then $\iota(P) \not\in \Supp(\fd)$ and

\item 
if $P \in \Supp(\fd)$ and $P = \iota(P)$, then $n_P = 1$.

\end{itemize}

\begin{defin}\label{def:semi-reduced}{\rm  Under the notation given above, we have the following definitions:
\begin{enumerate}
\item[(i)] An affine effective divisor $\fd$ on $C$ is said to be {\it semi-reduced} if $\fd = \fd_{\sr}$.
 
 \item[(ii)] A semi-reduced divisor $\fd$ is said to be {\it  $h$-reduced} if $\deg\fd \le g$.
 \end{enumerate}
 }
 \end{defin}
 
 \begin{rem}\label{rem:reduced}{\rm In \cite{galbraith, MWZ},  an $h$-reduced divisor defined as above
 is simply called a reduced divisor. On the other hand, the word ^^ reduced' is used  in a 
 different meaning in standard textbooks of algebraic geometry (e.g., \cite{iitaka}).
 Hence we use the word ^^ $h$-reduced' in order to avoid confusion.
 }
 \end{rem}

 We denote the  coordinate ring $K[x, y]/\langle y^2 - f \rangle$ (resp. $\Kbar[x, y]/\langle y^2 -f \rangle$)
 of the affine part of $C$ over $K$ by $K[C]$ (resp. over $\Kbar$ by $\Kbar[C]$).  For $g \in \Kbar[x, y]$,  we denote its class in 
 $\Kbar[C]$ by $[g]$.  For $P \in C$, $\mcO_P$ denotes the local ring at $P$ and $\ord_P$ means
 the discrete valuation at $P$.
 
 Let $\fd = \sum_{i=1}^r n_i P_i$, $P_i= (x_i, y_i)$ be a semi-reduced divisor on $C$. The divisor
 $\fd$ can be represented by a pair of  polynomials in $\Kbar[x]$ which is called the Mumford representation.
 We briefly explain it. We first recall the following lemma:
 
 \begin{lem}\label{lem:mumford-rep}{There exist unique polynomials $u(x), v(x) \in \Kbar[x]$  such that
  \begin{enumerate}
  \item[\rm{(i)}] $u(x) := \prod_{i=1}^r (x - x_i)^{n_i}$, 
  \item[\rm{(ii)}] $\deg v(x) < \deg u(x)$, $\ord_{P_i}([y - v(x)]) \ge n_i$, and  
  \item[\rm{(iii)}] $v(x)^2 - f$ is divisible by $u$.
  \end{enumerate}
  }
  \end{lem}
  For a proof,  see \cite[Lemma 10.3.5]{galbraith}.  
  
 Conversely if $u, v$ satisfying (i), (ii) and (ii) in Lemma~\ref{lem:mumford-rep} are given, 
 then we recover $\fd$ as the zero set  of the ideal generated by $[u]$ and 
 $[y - v]$
 in $\Kbar[C]$ counted with proper multiplicities  (see \cite{galbraith} in detail).

\begin{defin}\label{def:mumford-rep}{\rm For a semi-reduced divisor $\fd$, $(u, v)$ is
said to be the Mumford representation of $\fd$. Also, we denote the divisor given by a Mumford representation $(u,v)$ by $\fd(u,v)$. }
\end{defin}

\begin{rem}\label{rem:mumford-rep}{\rm
Define ideals $I(\fd)$ and  $\widetilde{I(\fd)}$ in $\Kbar[C]$ and $\Kbar[x, y]$, respectively 
as follows:
\begin{eqnarray*}
I(\fd) & := & \{\xi \in \Kbar[C] \mid \mbox{$\ord_P(\xi) \ge n_P$ for $\forall P \in \Supp(\fd)$}\}, \\
\widetilde{I(\fd)} &:= & \{ g \in \Kbar[x, y] \mid [g] \in \Kbar[C]\}.
\end{eqnarray*}
As it can be seen in \cite{takatoku20-2}, the Gr\"obner basis of $\widetilde{I(\fd)}$ with respect
to  the pure lexicographic oder $>_p$, $y>_p x$  is of the form
$\{u(x), y- v(x)\}$ where $(u, v)$ is the Mumford representation of $\fd$ and $\widetilde{I(\fd)} \cap  \Kbar[x] = \langle u(x) \rangle$. 
}
\end{rem}

\subsection{Semi-reduced divisors of degree $3$ on elliptic curves over $K$}

Let $E$ be an elliptic curve given by a Weierstrass equation
\[
y^2 = f(x) = x^3 + c_1x^2 + c_2x + c_3, \quad c_i \in K.
\]

If a semi-reduced divisor $\fd$ is defined over $K$, then the polynomials $u, v$ in the Mumford
representation of $\fd$ are in $K[x]$.
For later use, we consider the case of $\deg \fd = 3$. Let $\fd = P_1 + P_2 + P_3$ be a semi-reduced
divisor of degree $3$ ($P_i$'s are not necessarily distinct) and put $P_{\fd} := P_1\dot{+}P_2\dot{+}P_3$.

\begin{lem}\label{lem:deg3}{Assume that $P_{\fd} \neq O$ and let $(u_{\fd}, v_{\fd})$ be the Mumford representation of $\fd$. Then we have
\begin{enumerate}
\item[\rm (i)] $P_{\fd} \neq P_i$ ($i = 1, 2, 3$).
\item[\rm (ii)]  $\deg v_{\fd} = 2$.
\end{enumerate}
}
\end{lem}

\proof (i) By the definition of the addition on $E$, 
\[
P_1 + P_2 + P_3 - 3O \sim P_{\fd} - O.
\]
Hence, if $P_{\fd} = P_1$, then $P_2 + P_3 - 2O \sim 0$. This means $P_3 = \iota^*P_2$,  but this does not
occur as $\fd$ is semi-reduced. Similarly we have  $P_{\fd} \neq P_2, P_3$.

(ii) Since $P_{\fd} \neq O$, the points $P_1, P_2$ and $P_3$ are not  collinear. As the curve defined by $y - v_{\fd} = 0$ passes through
$P_1, P_2, P_3$ and $\deg v_{\fd} \le 2$, $\deg v_{\fd} = 2$.
\endproof

\begin{lem}\label{lem:deg3-2}{We keep the notation as before. Assume that $\fd$ is defined over $K$. Put $P_{\fd}:= (x_{\fd}, y_{\fd})$.
Then we have the following:
\begin{enumerate}
\item[\rm (i)] The point $P_{\fd}$ is a $K$-rational point of $E$, i.e., $x_{\fd}, y_{\fd} \in K$.
\item[\rm (ii)] The two polynomials $u_{\fd}, v_{\fd}$ satisfy  $u_{\fd}, v_{\fd} \in K[x]$. In particular, $v_{\fd}$ is of the form
\[
b_0(x - x_{\fd})(x - b_1) - y_{\fd} \quad (b_0, b_1 \in K).
\]
\end{enumerate}
}
\end{lem}

\proof (i) Since $P_{\fd}:= P_1\dot{+}P_2\dot{+}P_3$ and $\fd$ is defined over $K$, $P_{\fd}$ is
$\Gal(\Kbar/K)$-invariant. Hence $x_{\fd}, y_{\fd} \in K$.

(ii) The first statement follows from \cite[Lemma 10.3.10]{galbraith}.  We go on to show the second
statement.  By Lemma~\ref{lem:mumford-rep}, the divisor of $[y - v_{\fd}]$ is
$P_1 + P_2 + P_3 + P_4 - 4O$ for some $P_4 \in E$. Hence
\[
P_1 + P_2 + P_3 + P_4 +\iota^*P_4 - 5O \sim \iota^*P_4 - O.
\]
As $P_1 + P_2 + P_3 - 3O \sim P_{\fd} - O$ and $P_4 + \iota^*P_4 - 2O \sim 0$, 
$P_{\fd} - O \sim \iota^*P_4 - O$. Hence $\iota^*P_4 = P_{\fd}$, i.e., $P_4 = \iota^*P_{\fd} = (x_{\fd}, - y_{\fd})$
and $v_{\fd}$ has the desired form.
\endproof

\subsection{Trisections of elliptic surfaces and their Mumford representations}\label{subsec:2-2}
In what follows, we only consider the case when $C = \PP^1$, i.e., $S$ is an elliptic surface over $\PP^1$.

In \cite{bannai-tokunaga15, bannai-tokunaga17}, we made use of properiies of bisections of 
certain elliptic surfaces in order to construct Zarisiki $N$-ples. Our description for bisections there
 is nothing but the Mumford representation for bisections. One of our purposes of this article  is to consider a
 geometric application of Mumford representations of trisections along the line of \cite{bannai-tokunaga15,bannai-tokunaga17}.
 To this purpose, in this section, we consider the Mumford representations of trisections.
 Let $D$ be a trisection and $\fd_D$ be a divisor of degree $3$ obtained by the restriction of $D$ to
 $E_S$.  In this article we only consider trisections such that 
\begin{center}
   $\fd_D$ is a semi-reduced divisor on $E_S$.
 \end{center}
 We denote the rational point determined by $\fd_{D}$ by $P_D$, i.e, the rational point obtained by
 Theorerem~\ref{thm:shioda}. 
 In order to give the Mumford representation of $\fd_D$, we recall our setting in
  \cite[Sections 2.2.2, 2.2.3]{bannai-tokunaga15}. 
  Let $\Sigma_d$ ($d = \chi(\mcO_S)$) be
  the Hirzebruch surface of degree $d$. Take affine open subsets $U_1$ and $U_2$ of $\Sigma_d$ so that
  \begin{enumerate}
  \item[(i)] $U_i = \CC^2$ $(i = 1, 2)$ and 
  \item[(ii)] the coordinates of $U_1$ and $U_2$ are $(t, x)$ and $(s, x')$, respectively with $s = 1/t$ and
  $x' = x/t^d$.
  \end{enumerate}
  There exists a unique section $\Delta_0$ with respect to the ruling $\Sigma_d \to \PP^1$ such
  that $\Delta_0^2 = -d$. With the coordinates as abve, $\Delta_0$ is given by $x = x' = \infty$. 
  We also have a section given by $x = c_0t^d + c_1t^{d-1} + \ldots + c_d, c_i \in \CC$ on $U_1$ which is linearly equivalent to $\Delta_0 + d{\mathfrak f}$, $\mathfrak f$ being a fiber of the ruling $\Sigma_d \to \PP^1$.
  
  Under these settings, $S$ is obtained in the following way:
  
  \begin{enumerate}
  
  \item[(i)] There exists a reduced divisor $\mcT_S$ on $\Sigma_d$ such that (a) on $U_1$, $\mcT$ is given by 
  \[
  \mcT_S: f_{\mcT_S}(t, x) = x^3 + a_d(t)x^2 + a_{2d}(t)x + a_{3d}(t) = 0,
  \]
  where $a_{di}(t) \in \CC[t], \deg a_{di} \le di$ ($i = 1, 2, 3$), (b) $\mcT_S$ has at most simple singularities
  (see \cite{bpv} for simple singularities) and  (c) $E_S$ is given by the Weierstrass equation
  $y^2 = f_{\mcT}(t, x)$.
    
  \item[(ii)] The affine surface given by $y^2 = f_{\mcT_S}(t, x)$ can be extended to a double cover
  $f'_S:S' \to \Sigma_d$ and $S$ is obtained as the minimal resolution of $S'$. We denote the resolution
  by $\mu : S \to S'$ and put $f_S = f'_S \circ \mu$.  Note that $f_S(O) = \Delta_0$.
  See \cite{bannai-tokunaga15, tokunaga14, miranda-basic}, for detail. With the coordinates as above,
  $\CC(\PP^1) = \CC(t)$.
  
  \end{enumerate}
  
  \begin{rem}{\rm For the case of $d=2$, $S$ is a rational elliptic surface. If $S$ has a reducible singular
  fiber, $S$ can be obtained as $S_{\mcQ, z_o}$ for some quartic $\mcQ$ and a smooth point $z_o$ on
  $\mcQ$. If we choose coordinate of $\PP^2$ suitably, $f_{\mcT_S}(t, x)$ is a defining equation of $\mcQ$.
  See \cite[Section~2.2.2]{bannai-tokunaga15} for details.
  }
  \end{rem}

  We now consider the Mumford representation of $\fd_D$ via the Weiestrass equation of $E_S$ as above.
  Let $(u_D, v_D)$ be the Mumford representation of $\fd_{D}$ and we denote
  $P_D = (x_D, y_D) \in E_S(\CC(t))$. Since $\fd_D$ is defined over $\CC(t)$, $u_D, v_D \in \CC(t)[x]$.
  By Lemma~\ref{lem:deg3} and~\ref{lem:deg3-2}, $v_D$ is of the form
  \[
  b_0(x - x_D)(x - b_1) - y_D, \quad b_0, b_1 \in \CC(t),
  \]
  and
  \[
  (v_D)^2 - f_{\mcT_S} = b_0^2u_D(x - x_D).
  \]
  
 Conversely, assume that $P_o=(x_o, y_o) \in E_S(\CC(t))$ is given. Then for any polynomial $v \in \CC(t)[x]$ of 
 the form $b_0(x - x_o)(x - b_1) - y_o$  $b_0, b_1 \in \CC(t)$, $v^2 - f_{\mcT_S}$ has a decomposition
 $b_0^2(x - x_o)u$ where $u \in \CC(t)[x]$ is monic and $\deg u = 3$. Hence, $u$ and $v$ gives a divisor $\fd(u, v)$ of degree $3$ on $E_S$ 
 such that $\fd(u, v) - 3O \sim P_o - O$.  As $u$ is determined by $P_o$ and $b_0(x - x_o)(x - b_1) - y_o$, we denote the
 trisection given by $\fd(u, v)$ by $D(P_o, b_0, b_1)$. Note that  $\overline{\psi}(D(P_o, b_0, b_1)) = P_o$ holds and
 $\mcC_{D(P_o, b_0, b_1)}$ is given by $u = 0$.
 
 In the next section, we make use of this construction for trisections for a given $P_o$ with $x_o, y_o \in \CC[t]$ and
 $v$ of the form $b_o(x - x_o)(x - b_1(t)) - y_o$,  $b_o \in \CC^{\times}, b_ 1 \in \CC[t]$.

\section{Examples}\label{sec:examples}


We apply our observation in Section~\ref{sec:mumford-representation} to 
consider explicit examples for
trisections given by smooth cubics that appear in Proposition~\ref{prop:trisec-2}. As a result we obtain explicit examples of the Zariski pairs in Section~\ref{subsec:3-3}. We keep the notation as before.

\subsection{The case $\mcQ=C_o+L_1+L_2$}\label{subsec:6-1} 
We recall Example 5.2 in \cite{tokunaga14}.
Let $[T, X, Z]$ be homogeneous coordinates of $\PP^2$ and 
let $(t, x) = (T/Z, X/Z)$ be affine coordinates for $\CC^2 = \PP^2\setminus\{Z = 0\}$. Consider $\mcQ = C_o + L_1 + L_2$ given
by
\[
C_o: x - t^2 = 0, \quad L_1: x + 3t + 2 = 0, \quad L_2: x - 3t + 2 = 0.
\]
In this case, $p_o = [0, -2, 1]$ and we label $C_o\cap(L_1\cup L_2)$ by
$p_1=[-1, 1, 1]$, $p_2 = [-2, 4, 1]$, $p_3=[1, 1,1]$ and $p_4=[2, 4,1]$ (Note that we label $L_1$ and $L_2$ in a 
different way from \cite[Example 5.2]{tokunaga14}). The lines  $L_{13}$ and 
$L_{14}$ are given by $x - 1 = 0$ and $x - t - 2 = 0$, respectively. 
By choosing  $[0,1, 0]$ as $z_o$, we may assume that  $E_{\mcQ, z_o}$ is 
given by 
\[
y^2 = f_{\mcT_{S_{\mcQ, z_o}}} = (x - t^2)(x+ 3t + 2)(x - 3t+ 2).
\]
Under these settings, $P_{13}$ and $P_{14}$ are 
\[
P_{13} = (1, -3(t+1)(t-1)), \quad P_{14} = (t+2, -2\sqrt{2}(t+1)(t-2)).
\]
For simplicity we put $b_0= 1$. We need to choose $b_1$ appropriately so that $u = 0$ gives a smooth cubic and we have
the following:

\begin{enumerate}

\item[(i)] For $D(P_{13}, 1, c)$, $c \in \CC$, 
\[
u=x^3 - 2(c+1)x^2 + (7t^2 + c^2 + 2c - 11)x -(6c-14)t^2 - (c-3)^2.
\]
The curve $\mcC_{D(P_{13}, 1, c)}$ is given by $u = 0$ and it is smooth for general $c$.

\item[(ii)] For $D(P_{14}, 1, b_1)$ $b_1 =-(2\sqrt{2} + 3)t+ c)$, $c \in \CC$, 
\[
u=x^3 + a_1x^2 + a_2x + a_3,
\]
where
\begin{eqnarray*}
a_1 & := & (4\sqrt{2} + 5)t - 2c -3, \\
a_2 & := & (12\sqrt{2} + 12)t^2 - (4\sqrt{2}c +4c + 12\sqrt{2} + 13)t + c^2 + 4c -8\sqrt{2} - 6, \\
a_3 & := &-(-6c + 36\sqrt{2} + 36)t^2 - (c^2 - 12\sqrt{2}c -12c + 24\sqrt{2} + 40)t \\
      &  & -2c^2 + 8\sqrt{2}c - 16.
\end{eqnarray*}
The curve $\mcC_{D(P_{14}, 1, b_1)}$ is given by $u = 0$ and it is smooth for  general $c$.

\end{enumerate}
As $\displaystyle{[2]P_{14} = \left (\frac 98 t^2, -\frac{\sqrt{2}}{32}t(9t^2 - 16)\right )}$, $\mcC_{[2]P_{14}}$ is given by
$\displaystyle{x - \frac 98 t^2} = 0$. Put
\[
\mcB_1 = \mcQ + \mcC_{[2]P_{14}} +  \mcC_{D(P_{13}, 1, c)}, \quad
\mcB_2 = \mcQ + \mcC_{[2]P_{14}} + \mcC_{D(P_{14}, 1, b_1)}.
\]
Since $\mcB_1$ and $\mcB_2$ have the same combinatorics for general $c$, $(\mcB_1, \mcB_2)$ is a Zariski pair by Section~\ref{subsec:3-3}.

\subsection{The case $\mcQ=C_1+C_2$: The existence of curves satisfying   Theorem~\ref{thm:trisec-1}, (i) and (ii)}\label{subsec:6-2}

By choosing suitable coordinates,  we assume that $\mcQ$ is given by
$f_{\mcT_{S_{\mcQ, z_o}}} = 0$ and that the generic fiber $E_{\mcQ, z_o}$ is given by
$y^2 = f_{\mcT_{S_{\mcQ, z_o}}}$. 
 Let $P_{ij}$ and $Q_j$ be the points in considered in Section~\ref{subsec:case-2}. 
 Let $\fd_j$ ($j=1, 2, 3, 4$) be semi-reduced divisors
 \begin{eqnarray*}
 \fd_1 &:= &  P_{12} +P_{13} +P_{23},  \\
 \fd_2 &:= &  [-1]P_{12} + P_{13} +P_{23}, \\ 
 \fd_3 &:= & P_{12} + [-1]P_{13} +P_{23},  \\
 \fd_4 &:= & P_{12} + P_{13} + [-1]P_{23}.
 \end{eqnarray*}
 We give a observation about the Mumford representations of $\fd_j$ $(j = 1,2, 3,4)$.
 Let $(u_{\fd_j}, v_{\fd_j})$ be the Mumford representation of $\fd_j$. Put
 \[
 [\pm 1]P_{ij} = (x_{ij}(t), \pm y_{ij}(t)), \quad [\pm 1]Q_j = (x_{Q_j}(t), \pm y_{Q_j}(t)).
 \]
 As $C_{P_{ij}}$ and $C_{Q_j}$ are given by $x - x_{ij}(t) = 0$ and
 $x - x_{Q_j}(t) = 0$, respectively, $\deg x_{ij} = \deg x_{Q_j} = 1$.  Hence
 \[
 u_{\fd_j} = (x - x_{12}(t))(x - x_{13}(t))(x - x_{23}(t)), \quad  j = 1, 2, 3, 4, 
 \]
 i.e., $u_{\fd_j} = 0$ is $\triangle_{123}$ and $v_{\fd_j}$ satisfies
 \[
 (v_{\fd_j}^2 - f_{\mcT_{S_{\mcQ, z_o}}}) = \fd_j + [-1]Q_j - 4O,
 \]
 and $v_{\fd_j}^2 - f_{\mcT_{S_{Q, z_o}}} = c(x - x_{Q_j}(t))u_j, \, c \in {\CC(t)}^{\times}$.
 
 We now show the following Lemma which assures that the coefficient $c$ in the above equation is a constant.
 
 \begin{lem}\label{lem:mumford-rep2}{The polynomial 
 $v_{\fd_j}$ is of the form $v_{\fd_j}=c_0x^2 + c_1x + c_2$, $c_i \in \CC[t], \deg c_i = i$.
 In particular, $c$ in the above is a constant.
 }
 \end{lem}
 
 \proof  We prove the statement for $j=1$ only as the remaining cases can be proven
 in the same way. Put $\mcQ_{Q_1} = L_{12} + L_{13} + L_{23} + L_{Q_1}$ and
 $L_{Q_1}\cap \mcQ = \{q_1, q_2\}$. Consider a pencil $\Lambda_1$ of quartics given by
 $|\lambda \mcQ + \mu \mcQ_{Q_1}|_{[\lambda, \mu]\in \PP^1}$. $\Lambda_1$ satisfies the
 following properties:
 \begin{itemize}
 
 \item The base loucs of $\Lambda_1$ is $4(p_1 + p_2 + p_3) + 2(q_1 + q_2)$, where
 the coefficients of $p_i, q_j$'s mean the multiplicities.
 
 \item Since $\mcQ, \mcQ_{Q_1} \in \Lambda_1$, general members of  $\Lambda_1$ are irreducible.
 
 \item General members of $\Lambda_1$ have nodes at $p_i$ ($i = 1, 2, 3$).
 
 \end{itemize}
 
 Let $C_o$ be the unique conic passing through $p_i$ ($i = 1, 2, 3$) and $q_j$ ($j = 1, 2$).
 Choose $r \in C_o\setminus (\mcQ\cup\mcQ_{Q_1})$. Let $\mcQ_{\mu_r}$ be a member of $\Lambda_1$
 given by $\mcQ + \mu_r\mcQ_{Q_1}$ satisfying $r \in \mcQ_{\mu_r}$. 
 Then since the divisor on $C_o$ cut out by $\mcQ_{\mu_r}$ is  $C_o|_{\mcQ_{\mu_r}} =
 2(p_1 + p_2 + p_3) + q_1 + q_2 + r$, $C_o$ is an irreducible component of $\mcQ_{\mu_r}$. Write
 $\mcQ_{\mu_r} = C_o + C'_o$. Then $C'_o$ is irreducible and 
 $C'_o|_{\mcQ_{\mu_r}} =
 2(p_1 + p_2 + p_3) + q_1 + q_2$, i.e., $C_o^\prime = C_o$. Thus $2C_o \in \Lambda_1$. As $z_o \not\in C_o$, 
 $C_o$ is given by an equation of the form
 \[
 c_0x^2 + c_1x + c_2, \quad c_i \in \CC[t], \deg c_i = i.
 \]
 and, as $2C_o = \mcQ + \mu_r\mcQ_{Q_1}$, we have 
 \[
 (c_0 x^2 + c_1x + c_2)^2 = f_{\mcT_{S_{\mcQ, z_o}}} + \mu_r u_{\fd_1}( x - x_{Q_1}(t)).
 \]
 By the uniqueness of the Mumford representation, we have $v_{\fd_1} = c_0x^2 + c_1x + c_2$ and $c = \mu_r$. 
 \endproof
 

  By Lemma~\ref{lem:mumford-rep2} and $-y_{Q_j} = c_0x_{Q_j}^2 + c_1x_{Q_j} + c_2$,  we may assume
  $v_{\fd_j}$ is of the form
  \[
  b_0(x - x_{Q_j})(x + b_{10}t + b_{11}) - y_{Q_j}, \quad b_0 \in \CC^{\times},\, b_{1j} \in \CC,  j= 0, 1.
  \]
 For each $Q_j$, consider the one parameter family of polynomials over  $s \in \Delta_{\epsilon}, b_0(s) \neq 0$ given by
 \[
 v_{j, s} := b_0(s) (x - x_{Q_j})(x + b_{10}(s)t + b_{11}(s)) - y_{Q_j},
 \]
such that $v_{j, 0} = v_{\fd_j}$.  Since $v_{j, s}(x_{Q_j}(t)) = - y_{Q_j}$, we have
 \[
 v_{j, s}^2 - f_{\mcT_{S_{\mcQ, z_o}}} = b_0(s)^2(x - x_{Q_j})u_{j, s}(x).
 \]
 Hence $(u_{j, s}, v_{j, s})$ is a one parameter family of Mumford representations such that
 $\fd(u_{j,s}, v_{j, s}) - 3O \sim Q_j - O$. We also obtain a family of plane curves given by
 $u_{j, s}= 0$ for each $j$, which can be candidates of cubics in Theorem~\ref{thm:trisec-1}.
 We now go on to give an explicit example, which shows the existence of families of plane
 curves stated in Theorem~\ref{thm:trisec-1}.
 
 \begin{ex}\label{eg:2conics}{We use the same coordinates as in the case (i). Let $\mcQ  = C_1 + C_2$ given
 by 
 \[
 C_1 : x - t^2 = 0, \quad C_2: x^2 - 10tx + 25x - 36 = 0.
 \]
 In this case, $p_1=[3, 9, 1], p_2=[2, 4,1], p_3 =[6, 36, 1], p_4= [-1,1,1]$ and
 \[
 L_{12}: x - 5t + 6 = 0, \quad L_{13}: x - 9t + 18, \quad L_{23}: x - 8t + 12 = 0.
 \]
 By choosing  $[0,1, 0]$ as $z_o$, we may assume that  $E_{\mcQ, z_o}$ is 
given by 
\[
y^2 = f_{\mcT_{S_{\mcQ, z_o}}} = (x - t^2)(x^2 - 10tx + 25x - 36).
\]
 Under these settings, we have
 \[
 P_{12} =(5t - 6, -5(t-2)(t-3)),  P_{13}= (9t-18, -3(t-3)(t-6)), P_{23}=(8t-12, -4(t-2)(t-6)).
 \]
 Put $Q_1:= P_{12}\dot{+}P_{13}\dot{+}P_{23}, \, Q_2:=  [-1]P_{12}\dot{+}P_{13}\dot{+}P_{23}$. 
 Straightforward computations show
 \[
 Q_1 = (0, -6t), \quad Q_2= (10t-25, -6(t-5)).
 \]
 Hence $\mcQ$ has bitangents
 \[
 \mcL_1 : x = 0, \quad \mcL_2: x - 10t + 25 = 0.
 \]
 Put
 \begin{eqnarray*}
 v_{1, s} &:= & \frac 16 x\left (x - (11t - 36+s)\right )+ 6t \\
 v_{2, s} &:= & (x - 10t + 25)\left (x - (6t - 6+s)\right ) + 6(t-5),
 \end{eqnarray*}
where $s \in \CC$.
 By computing $v_{j, s}^2 - f_{\mcT_{S_{\mcQ, z_o}}}$ $(j = 1, 2)$, we have
 \[
 v_{1, s}^2 - f_{\mcT_{S_{\mcQ, z_o}}} = \frac 1{36}xu_{1, s}, \quad v_{2, s}^2 - f_{\mcT_{S_{\mcQ, z_o}}} = (x - 10t + 25)u_{2,s},
 \]
 where 
 \begin{eqnarray*}
 u_{1,s} & := & x^3 + (-2s - 22t + 36)x^2 + (s^2 + 22st + 157t^2 - 72s - 504t + 396)x\\ 
             & &  - 360t^3 + 72st + 1692t^2 - 2592t + 1296 \\
 u_{2,s} & := & x^3 - (2s + 22t - 36)x^2 - (-s^2 - 32st - 157t^2 + 62s + 504t - 396)x \\
   &  & - 10s^2t - 120st^2 - 360t^3 + 25s^2 + 432st + 1692t^2 - 360s - 2592t + 1296.
 \end{eqnarray*}
 Thus we have two families of Mumford representations $(u_{j,s}, v_{j,s})$ $(j = 1, 2)$ such that
 \begin{enumerate}
 
 \item[(i)] the divisor $\fd(u_{j,0}, v_{j, 0})$ satisfies $\fd(u_{j,0}, v_{j, 0}) = P_{12} + P_{13} + P_{23}$ and
 \item[(ii)] we have linear equivalences $\fd(u_{j, s}, v_{j, s}) - 3O \sim Q_j - O$  for $j = 1, 2$.
 
 \end{enumerate}
 
 Let $D_{1,s}:= D(Q_1, \frac 16, 11t - 36 + s)$ and $D_{2,s}:= D(Q_2, 1, 6t-6 + s)$ be trisections determined by
$ \fd(u_{1, s}, v_{1, s})$ and $\fd(u_{2, s}, v_{2, s})$, respectively.  
 Put  $E_{j, s} = \mcC_{D_{j,s}}$ $(j = 1, 2)$ be cubics given by $u_{j, s} = 0$ ($j = 1,2$), respectively. 
 Then we can check
 that  there exists   $\epsilon>0$, such that for $s\in \Delta_{\epsilon}$, $s\not=0$ the following hold: 
 \begin{enumerate}
 
 \item[(a)] The curves $E_{j, s}$ $(j = 1, 2)$ are smooth.
 
 \item[(b)] The curves $E_{j,s}$ $(j = 1,2)$ meet $\mcL_1$ at three distinct points.
 
 \item[(c)] The curves $E_{j,s}$ $(j = 1, 2)$ are tangent to $\mcQ$ at six distinct points.
 
 \end{enumerate}

 Let $\mcB_{j, s} := \mcQ + E_{j, s} + \mcL_1$ $(j = 1, 2)$. From the above facts and Section~3.3,  we see that the following statements hold:
 
 \begin{enumerate}
 
 \item[(i)] The curve $\mcB_{1, 0}$ is equal to $\mcB_{2,0}$.
 
 \item[(ii)] For each $s\not=0$,  $(\mcB_{1,s}, \mcB_{2,s})$ is a Zariski pair. 
 
 \end{enumerate}
 
 
 }
 \end{ex}


\bibliographystyle{abbrv}
\bibliography{biblio}

\noindent Shinzo BANNAI \\
National Institute of Technology, Ibaraki College\\
866 Nakane, Hitachinaka-shi, Ibaraki-Ken 312-8508 JAPAN \\
{\tt sbannai@gm.ibaraki-ct.ac.jp}\\

\noindent Ryosuke MASUYA, Non KAWANA and Hiro-o TOKUNAGA\\
Department of Mathematical  Sciences, Graduate School of Science, \\
Tokyo Metropolitan University, 1-1 Minami-Ohsawa, Hachiohji 192-0397 JAPAN \\
{\tt tokunaga@tmu.ac.jp}

\end{document}